\documentclass[12pt,twoside]{article}
\usepackage{amssymb}
\font\teneufm=eufm10 scaled \magstep1
\font\seveneufm=eufm7 scaled \magstep1
\font\fiveeufm=eufm5  scaled \magstep1
\newfam\eufmfam
\textfont\eufmfam=\teneufm
\scriptfont\eufmfam=\seveneufm
\scriptscriptfont\eufmfam=\fiveeufm

\usepackage{mathrsfs}

\newfam\msbfam
\font\tenmsb=msbm10 scaled \magstep1  \textfont\msbfam=\tenmsb
\font\sevenmsb=msbm7 scaled \magstep1 \scriptfont\msbfam=\sevenmsb
\font\fivemsb=msbm5 scaled \magstep1  \scriptscriptfont\msbfam=\fivemsb

\def\CC{{\mathbb C}}

\def\PP{{\mathbb P}}

\def\ra{\rightarrow}

 \def\HollowBoxx #1#2#3{{\dimen0=#1 \advance\dimen0 by -#2
       \dimen1=#1 \advance\dimen1 by #3
        \vrule height 0pt depth #3 width #2
       \hskip -#3
       \vrule height #1 depth #3 width #3}}
 \def\LeftContraction{\mathord{\kern1.45pt \HollowBoxx{6pt}{3.5pt}{.4pt}}\,}

 \def\HollowBox #1#2#3{{\dimen0=#1 \advance\dimen0 by -#3
       \dimen1=#1 \advance\dimen1 by #3
        \vrule height #1 depth #3 width #3
        \vrule height 0pt depth #3 width #2
        \hskip -#3}}
 \def\RightContraction{\mathord{\, \HollowBox{6pt}{3.1pt}{.4pt}} \kern1.6pt}

\def\qed{{\hfill $\Box$}}
\newtheorem{theorem}{THEOREM}[section]

\newtheorem{example}[theorem]{Example}
\newtheorem{remark}[theorem]{Remark}

\newtheorem{definition}[theorem]{Definition}

\begin{document}

\begin{center}
{\Large \bf Continuation of CR-Automorphisms 
\vspace{0.3cm}\\
of Levi Degenerate Hyperquadrics
\vspace{0.5cm}\\
to the Projective Space}
\vspace{0.5cm}\\
\normalsize A. V. Isaev and I. G. Kossovskiy
\end{center}

\begin{quotation} \small \sl 

\end{quotation}

\thispagestyle{empty}

\pagestyle{myheadings}
\markboth{A. V. Isaev and I. G. Kossovskiy}{CR-Automorphisms of Levi Degenerate Hyperquadrics}

\begin{quotation} \small\noindent\sl We show that every CR-automorphism of the closure of a Levi degenerate hyperquadric in the projective space extends to a holomorphic automorphism of the projective space.
\end{quotation}

\setcounter{section}{0}

\section{Introduction}
\setcounter{equation}{0}

Let $\langle z,z\rangle$ be a Hermitian form on $\CC^n$, where $z:=(z_1,\dots,z_n)$, and let $Q_{\langle \cdot,\cdot\rangle}$ be the quadric in $\CC^{n+1}$ associated to $\langle z,z\rangle$ as follows:
$$
Q_{\langle \cdot,\cdot\rangle}:=\left\{(z,w)\in\CC^{n+1}:\displaystyle\hbox{Im}\,w=\langle z,z\rangle\right\}.
$$
Consider the closure $\overline{Q}_{\langle \cdot,\cdot\rangle}$ of $Q_{\langle \cdot,\cdot\rangle}$ in $\CC\PP^{n+1}$. Clearly, we have
$$
\overline{Q}_{\langle \cdot,\cdot\rangle}=\left\{(\zeta:{\bf z}:{\bf w})\in\CC\PP^{n+1}:\displaystyle\frac{1}{2i}\left({\bf w}\overline\zeta-\zeta\overline{{\bf w}}\right)=\langle {\bf z},{\bf z}\rangle\right\},
$$
where $(\zeta:{\bf z}:{\bf w})$ are homogeneous coordinates in $\CC\PP^{n+1}$, the space $\CC^{n+1}$ is given in $\CC\PP^{n+1}$ by $\zeta\ne 0$ with $z={\bf z}/\zeta$, $w={\bf w}/\zeta$, and the hyperplane at infinity by $\zeta=0$. If the form $\langle z,z\rangle$ is non-degenerate, then $\overline{Q}_{\langle \cdot,\cdot\rangle}$ is a non-singular Levi non-degenerate hypersurface, and it is well-known that the following continuation phenomena hold for CR-automorphisms of $\overline{Q}_{\langle \cdot,\cdot\rangle}$: (i) every local $C^1$-smooth CR-automorphism of $\overline{Q}_{\langle \cdot,\cdot\rangle}$ extends to a global CR-automorphism of $\overline{Q}_{\langle \cdot,\cdot\rangle}$ [continuation along the hypersurface]; (ii) every global CR-automorphism of $\overline{Q}_{\langle \cdot,\cdot\rangle}$ extends to an automorphism of $\CC\PP^{n+1}$ [continuation away from the hypersurface] (see \cite{Po}, \cite{A}, \cite{Tan}). We note that continuation phenomena of both kinds for CR-automorphisms and, more generally, CR-isomorphisms of CR-manifolds other than $\overline{Q}_{\langle \cdot,\cdot\rangle}$ were observed by many authors (see e.g. \cite{Pi}, \cite{VEK}, \cite{NS}, \cite{Ka}, \cite{LS} and references therein).  

Let us now suppose that $\langle z,z\rangle$ is degenerate. Assuming without loss of generality that $\langle z,z\rangle$ is written as
$$
\langle z,z\rangle=\sum_{j=1}^m|z_j|^2-\sum_{j=m+1}^k|z_j|^2,
$$
where $m\ge 0$, $m\le k\le 2m$ and $k<n$, we see that the singular set $S_{\langle \cdot,\cdot\rangle}$ of $\overline{Q}_{\langle \cdot,\cdot\rangle}$ is
$$
S_{\langle \cdot,\cdot\rangle}=\left\{(0:{\bf z}:0)\in\CC\PP^{n+1}:{\bf z}_1=\dots={\bf z}_k=0\right\},
$$
and that the regular part $\tilde Q_{\langle \cdot,\cdot\rangle}:=\overline{Q}_{\langle \cdot,\cdot\rangle}\setminus S_{\langle \cdot,\cdot\rangle}$ of $\overline{Q}_{\langle \cdot,\cdot\rangle}$ is everywhere Levi degenerate. The collection of local CR-automorphisms of $\tilde Q_{\langle \cdot,\cdot\rangle}$ includes all maps of the form
$$
\begin{array}{llll}
z_j&\mapsto& z_j, & j=1,\dots,k,\\
\vspace{-0.3cm}\\
z_j&\mapsto& f_j(z), & j=k+1,\dots,n,\\
\vspace{-0.3cm}\\
w &\mapsto & w,
\end{array}
$$
where $f_j$ for $j=k+1,\dots,n$ are any functions holomorphic near the origin in $\CC^n$, such that the matrix $\left(\partial f_j/\partial z_l(0)\right)_{j,l=k+1,\dots,n}$ is non-degenerate. Thus, the first continuation phenomenon does not have any reasonable analogue in the Levi degenerate case. 

Interestingly, as we will see below, the second continuation phenomenon still holds in this case, provided one uses the right definition of (global) CR-automorphism of $\overline{Q}_{\langle \cdot,\cdot\rangle}$. Let $f: \overline{Q}_{\langle \cdot,\cdot\rangle}\ra\overline{Q}_{\langle \cdot,\cdot\rangle}$ be a bijective map. For $f$ to be called a CR-automorphism of $\overline{Q}_{\langle \cdot,\cdot\rangle}$ it is natural to require, first of all, that $f$ preserves the regular part $\tilde Q_{\langle \cdot,\cdot\rangle}$ and is a CR-automorphism of $\tilde Q_{\langle \cdot,\cdot\rangle}$ in the usual sense. It will turn out that this condition is sufficient for the second continuation phenomenon to hold if $k>m$, that is, if the Levi form of $\tilde Q_{\langle \cdot,\cdot\rangle}$ has eigenvalues of opposite signs at every point. Also, it is not hard to observe (see e.g. Example \ref{ex1} below) that for the case $k=m$ a condition on the behavior of $f$ on the singular set $S_{\langle \cdot,\cdot\rangle}$ is necessary. For instance, one can start by forcing $f$ to be continuous at the points of $S_{\langle \cdot,\cdot\rangle}$ (cf. the definition of CR-function on a singular quadratic cone given in \cite{CS}). However, as the following example shows, the continuity of $f$ on all of $\overline{Q}_{\langle \cdot,\cdot\rangle}$ is not sufficient for $f$ to extend to an automorphism of $\CC\PP^{n+1}$.

\begin{example}\label{ex1}\rm Let $k=m$. There exists an automorphism $\Phi$ of $\CC\PP^{n+1}$ that transforms $\overline{Q}_{\langle \cdot,\cdot\rangle}$ into the hypersurface  
$$
\overline{{\mathcal Q}}_{\langle \cdot,\cdot\rangle}:=\left\{(\zeta:{\bf z}:{\bf w})\in\CC\PP^{n+1}:|\zeta|^2-|{\bf w}|^2=\langle {\bf z},{\bf z}\rangle\right\}.
$$
Clearly, $\overline{{\mathcal Q}}_{\langle \cdot,\cdot\rangle}$ is the closure in $\CC\PP^{n+1}$ of the non-singular hypersurface 
$$
{\mathcal Q}_{\langle \cdot,\cdot\rangle}:=\left\{(z,w)\in\CC^{n+1}:\langle z,z\rangle+|w|^2=1\right\},
$$
which is the product of $\CC^{n-m}$ and the unit sphere in $\CC^{m+1}$. Under the map $\Phi$ the regular part $\tilde Q_{\langle \cdot,\cdot\rangle}$ is transformed into ${\mathcal Q}_{\langle \cdot,\cdot\rangle}$ and the singular part $S_{\langle \cdot,\cdot\rangle}$ into 
$$
{\mathcal S}_{\langle \cdot,\cdot\rangle}:=\overline{{\mathcal Q}}_{\langle \cdot,\cdot\rangle}\cap\{\zeta=0\}.
$$

We now define a map $f:\overline{{\mathcal Q}}_{\langle \cdot,\cdot\rangle}\ra\overline{{\mathcal Q}}_{\langle \cdot,\cdot\rangle}$ as follows: $f$ fixes every point of ${\mathcal S}_{\langle \cdot,\cdot\rangle}$, and on the finite part ${\mathcal Q}_{\langle \cdot,\cdot\rangle}$ it is the restriction of the automorphism of $\CC^{n+1}$ given by the formulas
$$
\begin{array}{llll}
z_j&\mapsto&z_j,&j=1,\dots,m,\\
\vspace{-0.3cm}\\
z_j&\mapsto&e^w z_j,&j=m+1,\dots,n,\\
\vspace{-0.3cm}\\
w&\mapsto&w.
\end{array}
$$
For $m+1\le l\le n$ let $U_l$ be the subset of $\CC\PP^{n+1}$ where ${\bf z}_l\ne 0$, and let $\rho^l:=\zeta/{\bf z}_l$, $\tau^l_j:={\bf z}_j/{\bf z}_l$, $\sigma^l:={\bf w}/{\bf z}_l$, with $j\ne l$, be coordinates in $U_l$. In these coordinates, on the intersection ${\mathcal Q}_{\langle \cdot,\cdot\rangle}\cap U_l$ the map $f$ is given by the formulas
$$
\begin{array}{lllll}
\rho^l&\mapsto&\displaystyle\rho^l e^{-\sigma^l/\rho^l},\\
\vspace{-0.3cm}\\
\tau^l_j&\mapsto&\displaystyle\tau^l_j e^{-\sigma^l/\rho^l},&j=1,\dots,m,\\
\vspace{-0.3cm}\\
\tau^l_j&\mapsto&\tau^l_j,&j=m+1,\dots,n,&j\ne l,\\
\vspace{-0.3cm}\\
\sigma^l&\mapsto&\displaystyle \sigma^l e^{-\sigma^l/\rho^l}.
\end{array}
$$
Letting $\rho^l\ra 0$ we see that $f$ is continuous on all of $\overline{{\mathcal Q}}_{\langle \cdot,\cdot\rangle}$, but does not extend holomorphically to a neighborhood of ${\mathcal S}_{\langle \cdot,\cdot\rangle}$.
\end{example}

The above example motivates imposing more restrictive conditions on the behavior of $f$ on the set $S_{\langle \cdot,\cdot\rangle}$. Our definition of CR-automorphism of $\overline{Q}_{\langle \cdot,\cdot\rangle}$ is therefore as follows.  

\begin{definition}\label{defcraut}\rm A map $f:\overline{Q}_{\langle \cdot,\cdot\rangle}\ra\overline{Q}_{\langle \cdot,\cdot\rangle}$ is called a CR-automorphism of $\overline{Q}_{\langle \cdot,\cdot\rangle}$ if $f$ is bijective and satisfies the conditions: (a) $f$ preserves $\tilde Q_{\langle \cdot,\cdot\rangle}$, and the restriction of $f$ to $\tilde Q_{\langle \cdot,\cdot\rangle}$ is a $C^1$-smooth CR-automorphism of $\tilde Q_{\langle \cdot,\cdot\rangle}$; (b) if $k=m$, then $f$ holomorphically extends to a neighborhood of $S_{\langle \cdot,\cdot\rangle}$.   
\end{definition}

In this short note the following result is obtained.

\begin{theorem}\label{main} \sl Every CR-automorphism of $\overline{Q}_{\langle \cdot,\cdot\rangle}$ extends to a holomorphic automorphism of $\CC\PP^{n+1}$.
\end{theorem}

The proof of Theorem \ref{main} given in Section \ref{proof} is inspired by some of the arguments utilized in \cite{NS}. These arguments can also be applied to yield the second continuation phenomenon for generic quadratic cones (see Remark \ref{cones}). We acknowledge that the proof for the case $k>m$ was suggested to us by S. Nemirovskii.

\section{Proof of Theorem \ref{main}}\label{proof}

Suppose first that $k>m$. Let $f$ be a CR-automorphism of $\overline{Q}_{\langle \cdot,\cdot\rangle}$. Since the Levi form of $\tilde Q_{\langle \cdot,\cdot\rangle}$ has eigenvalues of opposite signs at every point, $f$ extends to a biholomorphic map $f_1$ defined in a neighborhood $U$ of $\tilde Q_{\langle \cdot,\cdot\rangle}$. By a result of \cite{F}, \cite{Tak} (see also \cite{U}), the envelope of holomorphy of $U$ is either a Stein domain over $\CC\PP^{n+1}$ or coincides with $\CC\PP^{n+1}$. However, $\tilde Q_{\langle \cdot,\cdot\rangle}$ contains the copy of $\CC\PP^1$ given in $\CC\PP^{n+1}$ by the equations
$$
\begin{array}{l}
\zeta={\bf w},\\
\vspace{-0.3cm}\\
{\bf z}_1={\bf z}_k,\\
\vspace{-0.3cm}\\
{\bf z}_j=0,\quad\hbox{for}\,\, j=2,\dots,k-1,\,\,\hbox{and}\,\,j=k+1,\dots,n.
\end{array}
$$
Therefore, $U$ cannot lie in a Stein manifold, and thus the envelope of holomorphy of $U$ is $\CC\PP^{n+1}$.

Next, every locally biholomorphic map from an open subset $V$ of $\CC\PP^{n+1}$ into $\CC\PP^{n+1}$ extends to a locally biholomorphic map from the envelope of holomorphy of $V$ into $\CC\PP^{n+1}$ (see \cite{Ke}, \cite{I}, \cite{NS}). Therefore, $f_1$ extends to a locally biholomorphic map $f_2$ from $\CC\PP^{n+1}$ into $\CC\PP^{n+1}$. The compactness and simple connectedness of $\CC\PP^{n+1}$ now imply that $f_2$ is a holomorphic automorphism of $\CC\PP^{n+1}$, as required.

We now suppose that $k=m$ and use the notation introduced in Example \ref{ex1}. Choose an automorphism $\Phi$ as specified in this example, and let $g:=\Phi\circ f\circ\Phi^{-1}$. The restriction of $g$ to ${\mathcal Q}_{\langle \cdot,\cdot\rangle}$ is a $C^1$-smooth CR-automorphism of ${\mathcal Q}_{\langle \cdot,\cdot\rangle}$. 

Assume first that $m\ge 1$. In this case $g$ extends to a holomorphic automorphism of the domain     
$$
{\mathcal B}_{\langle \cdot,\cdot\rangle}:=\left\{(z,w)\in\CC^{n+1}:\langle z,z\rangle+|w|^2<1\right\},
$$
which is the product of $\CC^{n-m}$ and the unit ball in $\CC^{m+1}$. By assumption, $g$ also extends to a holomorphic map defined in a neighborhood of ${\mathcal S}_{\langle \cdot,\cdot\rangle}$. It then follows that $g$ extends to a biholomorphic map $g_1$ defined on the union $U$ of ${\mathcal B}_{\langle \cdot,\cdot\rangle}$ and a neighborhood of ${\mathcal S}_{\langle \cdot,\cdot\rangle}$. Observe that $U$ contains the copy of $\CC\PP^1$ given in $\CC\PP^{n+1}$ by the equations
$$
\begin{array}{l}
{\bf z}_j=0,\quad\hbox{for}\,\, j=1,\dots,n-1,\\
\vspace{-0.3cm}\\
{\bf w}=0.
\end{array}
$$
It follows that $U$ cannot lie in a Stein manifold, and thus by \cite{F}, \cite{Tak} the envelope of holomorphy of $U$ is $\CC\PP^{n+1}$. Arguing as earlier, we see that $g_1$ extends to a holomorphic automorphism of $\CC\PP^{n+1}$.

Suppose now that $m=0$. Clearly, the restriction of $g$ to ${\mathcal Q}_{\langle \cdot,\cdot\rangle}$ has the form
$$
w\mapsto\psi(w), \quad z\mapsto \Psi(z,w),
$$
where $\psi$ is a $C^1$-smooth diffeomorphism of the unit circle $S^1:=\{|w|=1\}$, and for every $w\in S^1$, the map $\Psi(z,w)$ is an automorphism of $\CC^n$. Each component of $\Psi(z,w)$ is given by a power series in $z$ whose coefficients are $C^1$-smooth functions on $S^1$. Since $g$ holomorphically extends to a neighborhood of ${\mathcal S}_{\langle \cdot,\cdot\rangle}$, for every $w\in S^1$ the map $\Psi(z,w)$ extends to an automorphism of $\CC\PP^n$. Hence $\Psi(z,w)$ is an affine map with respect to $z$ whose coefficients are $C^1$-smooth functions on $S^1$, i.e., for $j=1,\dots,n$ the $j$th component $\Psi_j(z,w)$ of $\Psi(z,w)$ has the form
$$
\Psi_j(z,w)=\sum_{l=1}^na_{jl}(w)z_j+b_j(w),
$$
where all the functions $a_{jl}(w)$, $b_j(w)$ are $C^1$-smooth on $S^1$.

Since $g$ holomorphically extends to a neighborhood of ${\mathcal S}_{\langle \cdot,\cdot\rangle}$, it is holomorphic on a set of the form
$$
\{(z,w)\in\CC^{n+1}: |z|>R, 1-\varepsilon<|w|<1+\varepsilon\},
$$
for some $R,\varepsilon>0$. Hence each of the functions $\psi(w), a_{jl}(w), b_j(w)$ holomorphically extends to a neighborhood of $S^1$, and therefore $g$ holomorphically extends to a neighborhood of ${\mathcal Q}_{\langle \cdot,\cdot\rangle}$ in $\CC^{n+1}$.

Thus $g$ continues to a biholomorphic map $g_1$ defined in a neighborhood $U$ of $\overline{{\mathcal Q}}_{\langle \cdot,\cdot\rangle}$. Since $\overline{{\mathcal Q}}_{\langle \cdot,\cdot\rangle}$ contains the copy of $\CC\PP^n$ given in $\CC\PP^{n+1}$ by the equation $\zeta={\bf w}$, it follows that $U$ cannot lie in a Stein manifold, and thus by \cite{F}, \cite{Tak} the envelope of holomorphy of $U$ is $\CC\PP^{n+1}$. Arguing as above, we see that $g_1$ extends to a holomorphic automorphism of $\CC\PP^{n+1}$.

The proof is complete.\qed

\begin{remark}\label{cones} \rm Let $\rho(z,w)$ be a real-valued homogeneous polynomial of degree two in $\hbox{Re}\,z_j$, $\hbox{Im}\,z_j$, $\hbox{Re}\,w$, $\hbox{Im}\,w$, $j=1,\dots,n$. The associated quadratic cone $C_{\rho}$ is defined to be the zero set of $\rho$. If $\rho$ is a Hermitian form on $\CC^{n+1}$, then the closure $\overline{C}_{\rho}$ of $C_{\rho}$ in $\CC\PP^{n+1}$ is projectively equivalent to a hypersurface $\overline{Q}_{\langle \cdot,\cdot\rangle}$ for which the corresponding form ${\langle z,z\rangle}$ is degenerate. If, on the other hand, $\rho$ is a generic non-Hermitian polynomial, then the singular set of $\overline{C}_{\rho}$ at infinity in $\CC\PP^{n+1}$ contains a compact complex variety. Therefore, the proof of Theorem \ref{main} given above yields the second continuation phenomenon for (appropriately defined) CR-automorphisms of $\overline{C}_{\rho}$ as well. For a study of CR-functions on irreducible quadratic cones and a complete classification of such cones in $\CC^2$ we refer the reader to \cite{CS} and references therein.
\end{remark}

{\obeylines
Department of Mathematics
The Australian National University
Canberra, ACT 0200
AUSTRALIA
E-mail: alexander.isaev@anu.edu.au, ilya.kossovskiy@anu.edu.au
}

\end{document}